\documentclass[english]{article}
\usepackage[T1]{fontenc}
\usepackage[latin9]{inputenc}
\usepackage{amsmath}

\makeatletter
\newcommand{\lyxaddress}[1]{
	\par {\raggedright #1
	\vspace{1.4em}
	\noindent\par}
}

\makeatother

\usepackage{babel}
\begin{document}
\title{Maximal generalization of Lanczos' derivative using one-dimensional
integrals}
\author{Andrej Liptaj\thanks{andrej.liptaj@savba.sk}}
\maketitle

\lyxaddress{\begin{center}
Institute of Physics, Slovak academy of Sciences\\
D\'{u}bravsk\'{a} cesta 9, 845 11 Bratislava, Slovakia
\par\end{center}}
\begin{abstract}
Derivative of a function can be expressed in terms of integration
over a small neighborhood of the point of differentiation, so-called
\emph{differentiation by integration} method. In this text a maximal
generalization of existing results which use one-dimensional integrals
is presented together with some interesting non-analytic weight functions.
\end{abstract}

\section{Introduction}

Cornelius Lanczos in his work \cite{lanczos} published a method of
differentiation by integration\footnote{The first person to publish such method was Cioranescu \cite{38_Cioranescu}.
The name of the method is however usually associated with Lanczos.}, where the derivative of a function is approximated by an integral.
The integral is performed over a small interval around the point of
differentiation with the approximation becoming exact in the limit
of the interval length approaching zero. For differentiable functions
one has
\[
f'\left(x_{0}\right)\equiv f'\left(x\right)|_{x=x_{0}}=\lim_{h\rightarrow0}\frac{3}{2h^{3}}\int_{-h}^{h}tf\left(x_{0}+t\right)dt.
\]
The expression is interesting from several aspects: it generalizes
the ordinary derivative\footnote{Converges in situations, where the ordinary derivative is not defined.}
and, also, its modifications might be useful for numerical differentiation
(see e.g. \cite{numDiff_1,numDiff_2}).

Since, the topic was addressed by several authors with noticeable
growth of interest in the last decade \cite{98_Gro,99_Shen,00_Hicks,05_Rangarajana,05_Burch,10_Wang,11_Gibaru,12_Diekema,18_Teruel}.
The millennial work \cite{00_Hicks} is probably the most interesting
of them: the authors actually provide a very broad generalization
of the Lanczos' formula for the first derivative and their approach
can be further and straightforwardly generalized to higher orders
(as done in this article). Their text is, surprisingly, widely overlook
by later works with exception of \cite{numDiff_2,05_Burch,12_Diekema},
which, however, do not exploit the potential of it.

In what follows, the second section will be dedicated to the generalization
of the Lanczos' approach for the first derivative. The next section
will cover generalization to higher-order derivatives and, in the
fourth section, a short discussion will follow. A summary and conclusion
will constitute the last section.

Let me remark that generalizations based on multidimensional integrals
can be found in literature (e.g. formula 2.31 in \cite{18_Teruel}).
Unlike other approaches, they, presumably, do not represent a special
case of the generalization presented here and remain an independent
way of generalizing the Lanczos' derivative.

\section{First derivative\label{sec:First-derivative}}

Let me restate the findings from \cite{00_Hicks}. The key observation
which allows for large generalizations is, that the approximation
of the derivative can be seen as averaging the derivative over some
small interval $\left[x_{0}-h,x_{0}+h\right]$ around the point of
differentiation $x_{0}$. This average might be understood as weighted
average with a weight function $w_{h}$
\begin{align}
f'\left(x_{0}\right) & \approx\int_{x_{0}-h}^{x_{0}+h}w_{h}\left(t\right)f'\left(t\right)dt,\label{eq:f1_h}\\
\text{where} & \int_{x-h}^{x+h}w_{h}\left(t\right)dt=1.\nonumber 
\end{align}
Negative weights cannot be excluded, yet the condition $0\leq w_{h}\left(t\right)$
might be adopted if desired. Because the weight functions are of the
most interest here, a modified version of (\ref{eq:f1_h}) will be
used throughout this text
\begin{align*}
f'\left(x_{0}\right) & \approx\int_{-1}^{1}w\left(t\right)f'\left(x_{0}+ht\right)dt,\\
 & \int_{-1}^{1}w\left(t\right)dt=1,
\end{align*}
so that the weight functions are defined on a ``standard'' interval
$\left[-1,1\right]$. One has
\[
w\left(t\right)=hw_{h}\left(x_{0}+ht\right).
\]
Using integration by partes one arrives to
\[
f'\left(x_{0}\right)\approx\frac{1}{h}\left[w\left(t\right)f\left(x_{0}+ht\right)\right]_{t=-1}^{t=1}-\frac{1}{h}\int_{-1}^{1}w'\left(t\right)f\left(x_{0}+ht\right)dt.
\]
Two interesting observations can be done:
\begin{itemize}
\item If $w$ is constant $w\left(t\right)=0.5$ then the standard definition
of the derivative is recovered
\[
f'\left(x_{0}\right)\approx\frac{f\left(x_{0}+h\right)-f\left(x_{0}-h\right)}{2h}.
\]
\item If $w\left(-1\right)=w\left(+1\right)=0$ then a \emph{differentiation
by integration} method is constructed
\[
f'\left(x_{0}\right)\approx-\frac{1}{h}\int_{-1}^{1}w'\left(t\right)f\left(x_{0}+ht\right)dt.
\]
\end{itemize}
The usual Lanczos' expression is obtained for 
\[
w\left(t\right)=\frac{3}{4}\left(1-t^{2}\right).
\]
Indeed
\[
w'\left(t\right)=-\frac{3}{2}t\quad\longrightarrow\quad f'\left(x_{0}\right)\approx\frac{3}{2h}\int_{-1}^{1}tf\left(x_{0}+ht\right)dt=\frac{3}{2h^{3}}\int_{-h}^{h}zf\left(x_{0}+z\right)dz.
\]
At this point one can formulate the generalization: \emph{Any} differentiable
function $w$ which satisfies
\[
\int_{-1}^{1}w\left(t\right)dt=1\text{ and }w\left(-1\right)=w\left(1\right)=0
\]
can be used for \emph{differentiation by integration} in the following
manner
\[
f'\left(x_{0}\right)\approx-\frac{1}{h}\int_{-1}^{1}w'\left(t\right)f\left(x_{0}+ht\right)dt,
\]
where its derivative $w'$ appears.

Let me define some useful terms: ``kernel function'' will from now
on refer to the function which is being integrated (together with
function values) in the \emph{differentiation by integration} procedure\footnote{In case of the first derivative the kernel function is $w'.$}
and let me note by small zero those anti-derivatives of a function
$k$ which take value zero at minus one
\[
k_{0}^{\left(-n\right)}\left(t\right)|_{t=-1}=0,\quad\frac{d}{dt}k_{0}^{\left(-n\right)}=k_{0}^{\left(-n+1\right)}.
\]
One can now address the question about a proper kernel function (inverse
implication). From what was show one can deduce: $k$ is valid kernel
function iff
\begin{equation}
k_{0}^{\left(-1\right)}\left(t\right)|_{t=+1}=0\text{ and }\int_{-1}^{1}k_{0}^{\left(-1\right)}\left(t\right)dt=1.\label{eq:Der1_iff}
\end{equation}
The first of the two conditions is equivalent to
\[
\int_{-1}^{1}k\left(t\right)dt=0.
\]
Indeed, proceeding by integration by parts one observes ($\lambda=\varLambda'$)
\[
-\frac{1}{h}\int_{-1}^{1}\lambda\left(t\right)f\left(x_{0}+ht\right)dt=-\frac{1}{h}\left[\varLambda\left(t\right)f\left(x_{0}+ht\right)\right]_{t=-1}^{t=1}+\int_{-1}^{1}\varLambda\left(t\right)f'\left(x_{0}+ht\right)dt.
\]
If $\varLambda\left(\pm1\right)\neq0$, one cannot make vanish the
first term on the RHS for a general function $f$. If one takes the
limit $h\rightarrow0$ in the second term (using continuity of $f'$)
one arrives to
\[
\lim_{h\rightarrow0}\int_{-1}^{1}\varLambda\left(t\right)f'\left(x_{0}+ht\right)dt=\int_{-1}^{1}\varLambda\left(t\right)f'\left(x_{0}\right)dt=f'\left(x_{0}\right)\int_{-1}^{1}\varLambda\left(t\right)dt.
\]
One sees that a function with integral different from one provides
wrong value of the derivative. Formulas (\ref{eq:Der1_iff}) express
sufficient and necessary conditions a kernel function has to fulfill,
they represent the largest possible generalization of the Lanczos'
approach.

\section{Higher order derivatives}

\subsection{Main result}

Repeated integration by parts allows for immediate generalization
\begin{align*}
 & \int_{-1}^{1}w\left(t\right)f^{\left(n\right)}\left(x_{0}+ht\right)dt=\\
 & =\frac{1}{h}\left[w\left(t\right)f^{\left(n-1\right)}\left(x_{0}+ht\right)\right]_{-1}^{1}-\frac{1}{h}\int_{-1}^{1}w'\left(t\right)f^{\left(n-1\right)}\left(x_{0}+ht\right)dt\\
 & =\frac{1}{h}\left[w\left(t\right)f^{\left(n-1\right)}\left(x_{0}+ht\right)\right]_{-1}^{1}-\frac{1}{h^{2}}\left[w'\left(t\right)f^{\left(n-2\right)}\left(x_{0}+ht\right)\right]_{-1}^{1}\\
 & \qquad+\frac{1}{h^{2}}\int_{-1}^{1}w''\left(t\right)f^{\left(n-2\right)}\left(x_{0}+ht\right)dt\\
 & =\left(\frac{-1}{h}\right)^{n}\int_{-1}^{1}w^{\left(n\right)}\left(t\right)f\left(x_{0}+ht\right)dt+\sum_{k=0}^{n-1}\frac{\left(-1\right)^{k}}{h^{k+1}}\left[w^{\left(k\right)}\left(t\right)f^{\left(n-1-k\right)}\left(x_{0}+ht\right)\right]_{-1}^{1}
\end{align*}
To make, for a general function $f$, the second term vanish, one
has to require
\begin{equation}
w^{\left(k\right)}\left(-1\right)=w^{\left(k\right)}\left(1\right)=0\text{ for all }k=0,1,\ldots,n-1.\label{eq:conditionZero}
\end{equation}
Having this property, then, with appropriate weight function
\begin{equation}
\int_{-1}^{1}w\left(t\right)dt=1,\label{eq:conditionArea}
\end{equation}
and assuming $f^{\left(n\right)}$ is continuous, one interprets the
first term as an approximation of the $n$-th derivative
\begin{align*}
\lim_{h\rightarrow0}\left(\frac{-1}{h}\right)^{n}\int_{-1}^{1}w^{\left(n\right)}\left(t\right)f\left(x_{0}+ht\right)dt & =\lim_{h\rightarrow0}\int_{-1}^{1}w\left(t\right)f^{\left(n\right)}\left(x_{0}+ht\right)dt\\
 & =f^{\left(n\right)}\left(x_{0}\right)\int_{-1}^{1}w\left(t\right)dt\\
 & =f^{\left(n\right)}\left(x_{0}\right).
\end{align*}
Like at the end of the Sec. \ref{sec:First-derivative}, one can inverse
the whole procedure, start with expression $\int_{-1}^{1}w^{\left(n\right)}\left(t\right)f\left(x_{0}+ht\right)dt$
and proceed to $n$ repeated integrations by parts (integrate $w^{\left(n\right)}$
and differentiate $f$). As result one can immediately conclude: If
$k$ is to be a valid kernel for \emph{differentiation by integration}
in the formula
\begin{equation}
f^{\left(n\right)}\left(x_{0}\right)\approx\left(\frac{-1}{h}\right)^{n}\int_{-1}^{1}k\left(t\right)f\left(x_{0}+ht\right)dt\label{eq:DerN_gen}
\end{equation}
then 
\[
k_{0}^{\left(-n\right)}\left(1\right)=0\text{ for all }n=1,2,\ldots,n
\]
and
\[
\int_{-1}^{1}k_{0}^{\left(-n\right)}\left(t\right)dt=1.
\]
With these statements valid for \emph{any} weight/kernel functions
for which appropriate derivatives/integrals exist, one can claim that,
for the Lanczos' derivative written in the from (\ref{eq:DerN_gen}),
the generalization is maximal.

\subsection{Examples}

With the acquired knowledge one can propose some new, potentially
interesting kernels and weight functions. Idea of universality might
be a compelling one, by which I mean the independence on the order
of the derivative (from now on noted $n$). Kernels have to be $n$-dependent\footnote{The LHS of (\ref{eq:DerN_gen}) is $n$-dependent, so has to be the
RHS. But, with the exception of $\left(\frac{-1}{h}\right)^{n}$,
there are no other explicit $n$-dependent factors on the RHS, thus
the dependence must be hidden in $k\left(t\right)$.}, but one can look for $n$-independent weight functions. Such a universal
weight function has to fulfill condition (\ref{eq:conditionZero})
for all derivatives, yet it cannot be zero so as to respect the condition
(\ref{eq:conditionArea}). Therefore it must be non-analytic at -1
and 1.

As first example I propose
\[
w_{e}=\frac{1}{K}\exp\left(\frac{1}{x^{2}-1}\right)\text{ with }K\approx0.4439938161680786.
\]
With no explicit $n$-dependence in the weight function, this dependence
comes from differentiation
\[
f^{\left(n\right)}\left(x_{0}\right)\approx\left(\frac{-1}{h}\right)^{n}\frac{1}{K}\int_{-1}^{1}dtf\left(x_{0}+ht\right)\frac{d^{n}}{dt^{n}}\exp\left(\frac{1}{\left[t^{2}-1\right]}\right)
\]
Explicit formulas for the first three derivatives are
\[
f^{'}\left(x_{0}\right)\approx\frac{2}{hK}\int_{-1}^{1}dtf\left(x_{0}+ht\right)\frac{t}{\left(t-1\right)^{2}\left(t+1\right)^{2}}\exp\left(\frac{1}{\left[t^{2}-1\right]}\right)
\]
\[
f^{''}\left(x_{0}\right)\approx\frac{2}{h^{2}K}\int_{-1}^{1}dtf\left(x_{0}+ht\right)\frac{3t^{4}-1}{\left(t-1\right)^{4}\left(t+1\right)^{4}}\exp\left(\frac{1}{\left[t^{2}-1\right]}\right)
\]
\[
f^{'''}\left(x_{0}\right)\approx\frac{4}{h^{3}K}\int_{-1}^{1}dtf\left(x_{0}+ht\right)\frac{x\left(6t^{6}+3t^{4}-10t^{2}+3\right)}{\left(t-1\right)^{6}\left(t+1\right)^{6}}\exp\left(\frac{1}{\left[t^{2}-1\right]}\right)
\]
Even more interesting example is a one with shifted Fabius function
\cite{66_Fabius}
\[
w_{Fb}\left(t\right)=Fb(t+1).
\]
 The Fabius function (which I note $Fb$) is non-analytic for all
$0\leq x$ and its behavior with respect to the conditions (\ref{eq:conditionZero},\ref{eq:conditionArea})
can be deduced from differential functional equation 
\begin{equation}
Fb'\left(x\right)=2Fb\left(2x\right).\label{eq:FabEq}
\end{equation}
One has
\[
\int_{0}^{2}Fb\left(x\right)dx\stackrel{x=2z}{=}\int_{0}^{1}2Fb\left(2z\right)dz=\int_{0}^{1}Fb'\left(z\right)dz=\left[Fb\left(z\right)\right]_{z=0}^{z=1}=1
\]
\[
0=Fb\left(0\right)=\frac{1}{2}Fb'\left(0\right)=\frac{1}{2}\frac{1}{4}Fb''\left(0\right)=\ldots
\]
\[
Fb^{\left(n\right)}\left(2\right)=\frac{1}{2^{n+1}}Fb^{\left(n+1\right)}\left(1\right)=0,
\]
where the very last equality (all derivatives vanishing at $x=1$)
is consequence of the symmetry condition $Fb\left(1-x\right)=1-Fb(x)$
and the behavior of derivatives at $x=0$. When shifting Fabius function
to the interval $\left[-1,1\right]$ all mentioned properties remain
conserved (on the shifted the interval). Equation (\ref{eq:FabEq})
allows us to formulate corresponding kernel functions in a very elegant
way, where the explicit dependence on derivatives is not present\footnote{One can notice that the expression is defined for any real value of
$n$.}
\[
f^{\left(n\right)}\left(x_{0}\right)\approx\left(\frac{-1}{h}\right)^{n}2^{\frac{1}{2}n\left(n+1\right)}\int_{-1}^{1}Fb\left[2^{n}\left(t+1\right)\right]f\left(x_{0}+ht\right)dt.
\]
Value of the Fabius function for $1<x$ can be very easily related
to the value of this function on the interval $\left[0,1\right]$.
Using an efficient method\footnote{Use of tabulated values, or recipes from \cite{Fb_Rec1,Fb_Rec2,Fb_Rec3}.}
for its evaluation on the interval $\left[0,1\right]$, one achieves
an effective method for computing kernel function values and thus
the whole integral, and this for any order of the derivative.

\section{Discussion}

One of the most cited results \cite{05_Rangarajana} generalizes the
Lanczos' derivative by using Legendre polynomials\footnote{A similar result was in the same year published by \cite{05_Burch}.}.
It might be interesting to check its behavior from the perspective
of presented results. The authors of \cite{05_Rangarajana} propose
(among others) the following form of the kernel function\footnote{Factor $\left(-1\right)^{n}$ is here to cancel the same factor in
(\ref{eq:DerN_gen}) from in front of the integral.}
\[
k_{n}\left(x\right)=\frac{\left(-1\right)^{n}}{2}\left(2n+1\right)!!P_{n}\left(x\right),
\]
with $P_{n}\left(x\right)$ being the Legendre polynomials. The latter
can be defined by Rodrigues' formula
\[
P_{n}\left(x\right)=\frac{1}{2^{n}n!}\frac{d^{n}}{dx^{n}}\left(x^{2}-1\right)^{n}.
\]
Observing the inner bracket (going to zero for $x=\pm1$) being raised
to the $n$-th power, one immediately sees that the condition (\ref{eq:conditionZero})
is obeyed. Next, one can study the integral of the weight function
\[
\frac{\left(-1\right)^{n}}{2}\left(2n+1\right)!!\frac{1}{2^{n}n!}\int_{-1}^{1}\left(x^{2}-1\right)^{n}dx.
\]
With partial results \cite{wolfram}
\[
\int_{-1}^{1}\left(x^{2}-1\right)^{n}dx=\sqrt{\pi}\left(-1\right)^{n}\frac{n!}{\varGamma\left(n+\frac{3}{2}\right)}\text{ and }\varGamma\left(n+\frac{3}{2}\right)=\sqrt{\pi}\frac{\left(2n+1\right)!!}{2^{\left(n+1\right)}}
\]
one finds that also the condition (\ref{eq:conditionArea}) is respected.

Several other realizations for \emph{differentiation by integration}
can be found in the literature, most of them with higher technical
complexity then the previous one. From what was shown, all of these
representation (based on one-dimensional integrals) have to comply
with the restrictions (\ref{eq:conditionZero}) and (\ref{eq:conditionArea}).

This text focuses on the main result of generalizing the Lanczos'
derivative and does not address specific issues of precision and rapidity
of convergence in case of a numerical implementation and related questions
of kernel function preference. With kernel function being completely
general (possibly non-analytic everywhere) one can hardly rely on
standard tools for error estimates (i.e. Taylor series). In any specific
context the recipes exiting in the literature are to be used.

\section{Summary, conclusion}

In this text the result published in \cite{00_Hicks} was generalized
to higher-order derivatives and, assuming pattern (\ref{eq:DerN_gen}),
this generalization is maximal. Restrictions (\ref{eq:conditionZero})
and (\ref{eq:conditionArea}) allow for a very broad family of functions,
which might make the search for well performing kernels for numerical
purposes more efficient.

\end{document}